\documentclass{amsart}
\usepackage{amsfonts}

\newtheorem{theorem}{Theorem}
\theoremstyle{plain}

\newtheorem{corollary}{Corollary}

\newtheorem{lemma}{Lemma}

\newtheorem{proposition}{Proposition}
\newtheorem{remark}{Remark}

\numberwithin{equation}{section}
\input{tcilatex}

\begin{document}
\title[Quadratic Reverses of the Continuous Triangle Inequality]{Quadratic
Reverses of the Continuous Triangle Inequality for Bochner Integral of
Vector-Valued Functions in Hilbert Spaces}
\author{Sever S. Dragomir}
\address{School of Computer Science and Mathematics\\
Victoria University of Technology\\
PO Box 14428, MCMC 8001\\
Victoria, Australia.}
\email{sever@csm.vu.edu.au}
\urladdr{http://rgmia.vu.edu.au/SSDragomirWeb.html}
\date{April 16, 2004.}

\begin{abstract}
Some quadratic reverses of the continuous triangle inequality for Bochner
integral of vector-valued functions in Hilbert spaces are given.
Applications for complex-valued functions are provided as well.
\end{abstract}

\maketitle

\section{Introduction}

Let $f:\left[ a,b\right] \rightarrow \mathbb{K}$, $\mathbb{K}=\mathbb{C}$ or 
$\mathbb{R}$ be a Lebesgue integrable function. The following inequality is
the continuous version of the triangle inequality%
\begin{equation}
\left\vert \int_{a}^{b}f\left( x\right) dx\right\vert \leq
\int_{a}^{b}\left\vert f\left( x\right) \right\vert dx  \label{1.1}
\end{equation}%
and plays a fundamental role in Mathematical Analysis and its applications.

It seems, see \cite[p. 492]{MPF}, that the first reverse inequality for (\ref%
{1.1}) was obtained by J. Karamata in his book from 1949, \cite{K}:%
\begin{equation}
\cos \theta \int_{a}^{b}\left\vert f\left( x\right) \right\vert dx\leq
\left\vert \int_{a}^{b}f\left( x\right) dx\right\vert   \label{1.2}
\end{equation}%
provided%
\begin{equation*}
\left\vert \arg f\left( x\right) \right\vert \leq \theta ,\ \ x\in \left[ a,b%
\right] ,
\end{equation*}%
where $\theta $ is a given angle in $\left( 0,\frac{\pi }{2}\right) .$

This integral inequality is the continuous version of a reverse inequality
for the generalised triangle inequality%
\begin{equation}
\cos \theta \sum_{i=1}^{n}\left\vert z_{i}\right\vert \leq \left\vert
\sum_{i=1}^{n}z_{i}\right\vert ,  \label{1.3}
\end{equation}%
provided%
\begin{equation*}
a-\theta \leq \arg \left( z_{i}\right) \leq a+\theta ,\ \ \text{for \ }i\in
\left\{ 1,\dots ,n\right\} ,
\end{equation*}%
where $a\in \mathbb{R}$ and $\theta \in \left( 0,\frac{\pi }{2}\right) ,$
which, as pointed out in \cite[p. 492]{MPF}, was first discovered by M.
Petrovich in 1917, \cite{P}, and, subsequently rediscovered by other
authors, including J. Karamata \cite[p. 300 -- 301]{K}, H.S. Wilf \cite{W},
and in an equivalent form by M. Marden \cite{M}.

The first to consider the problem for sums in the more general case of
Hilbert and Banach spaces, were J.B. Diaz and F.T. Metcalf \cite{DM}.

In our previous work \cite{SSD1}, we pointed out some continuous versions of
Diaz and Metcalf results providing reverses of the generalised triangle
inequality in Hilbert spaces.

We mention here some results from \cite{SSD1} which may be compared with the
new ones obtained in Sections 2 and 3 below.

\begin{theorem}
\label{ta}If $f\in L\left( \left[ a,b\right] ;H\right) ,$ the space of
Bochner integral functions defined on $\left[ a,b\right] $ and with values
in the Hilbert space $H,$ and there exists a constant $K\geq 1$ and a vector 
$e\in H,$ $\left\Vert e\right\Vert =1$ such that%
\begin{equation}
\left\Vert f\left( t\right) \right\Vert \leq K\func{Re}\left\langle f\left(
t\right) ,e\right\rangle \ \ \ \text{for a.e. }t\in \left[ a,b\right] ,
\label{1.4}
\end{equation}%
then we have the inequality:%
\begin{equation}
\int_{a}^{b}\left\Vert f\left( t\right) \right\Vert dt\leq K\left\Vert
\int_{a}^{b}f\left( t\right) dt\right\Vert .  \label{1.5}
\end{equation}%
The case of equality holds in (\ref{1.5}) if and only if%
\begin{equation}
\int_{a}^{b}f\left( t\right) dt=\frac{1}{K}\left( \int_{a}^{b}\left\Vert
f\left( t\right) \right\Vert dt\right) e.  \label{1.6}
\end{equation}
\end{theorem}

As particular cases of interest that may be applied in practice, we note the
following corollaries established in \cite{SSD1}.

\begin{corollary}
\label{cb}Let $e$ be a unit vector in the Hilbert space $\left(
H;\left\langle \cdot ,\cdot \right\rangle \right) ,$ $\rho \in \left(
0,1\right) $ and $f\in L\left( \left[ a,b\right] ;H\right) $ so that%
\begin{equation}
\left\Vert f\left( t\right) -e\right\Vert \leq \rho \ \ \ \text{for a.e. }%
t\in \left[ a,b\right] .  \label{1.7}
\end{equation}%
Then we have the inequality%
\begin{equation}
\sqrt{1-\rho ^{2}}\int_{a}^{b}\left\Vert f\left( t\right) \right\Vert dt\leq
\left\Vert \int_{a}^{b}f\left( t\right) dt\right\Vert ,  \label{1.8}
\end{equation}%
with equality if and only if%
\begin{equation}
\int_{a}^{b}f\left( t\right) dt=\sqrt{1-\rho ^{2}}\left(
\int_{a}^{b}\left\Vert f\left( t\right) \right\Vert dt\right) \cdot e.
\label{1.9}
\end{equation}
\end{corollary}

\begin{corollary}
\label{cc}Let $e$ be a unit vector in $H$ and $M\geq m>0.$ If $f\in L\left( %
\left[ a,b\right] ;H\right) $ is such that%
\begin{equation}
\func{Re}\left\langle Me-f\left( t\right) ,f\left( t\right) -me\right\rangle
\geq 0\ \   \label{1.10}
\end{equation}%
or, equivalently,%
\begin{equation}
\left\Vert f\left( t\right) -\frac{M+m}{2}e\right\Vert \leq \frac{1}{2}%
\left( M-m\right) \ \ \   \label{1.11}
\end{equation}%
for a.e. $t\in \left[ a,b\right] ,$ then we have the inequality%
\begin{equation}
\frac{2\sqrt{mM}}{M+m}\int_{a}^{b}\left\Vert f\left( t\right) \right\Vert
dt\leq \left\Vert \int_{a}^{b}f\left( t\right) dt\right\Vert ,  \label{1.12}
\end{equation}%
or, equivalently%
\begin{equation}
0\leq \int_{a}^{b}\left\Vert f\left( t\right) \right\Vert dt-\left\Vert
\int_{a}^{b}f\left( t\right) dt\right\Vert \leq \frac{\left( \sqrt{M}-\sqrt{m%
}\right) ^{2}}{M+m}\left\Vert \int_{a}^{b}f\left( t\right) dt\right\Vert .
\label{1.13}
\end{equation}%
The equality holds in (\ref{1.12}) (or in the second part of (\ref{1.13}))
if and only if%
\begin{equation*}
\int_{a}^{b}f\left( t\right) dt=\frac{2\sqrt{mM}}{M+m}\left(
\int_{a}^{b}\left\Vert f\left( t\right) \right\Vert dt\right) e.
\end{equation*}
\end{corollary}

The case of additive reverse inequalities for the continuous triangle
inequality has been considered in \cite{SSD2}.

We recall here the following general result.

\begin{theorem}
\label{td}If $f\in L\left( \left[ a,b\right] ;H\right) $ is such that there
exists a vector $e\in H,$ $\left\Vert e\right\Vert =1$ and $k:\left[ a,b%
\right] \rightarrow \lbrack 0,\infty )$ a Lebesgue integrable function such
that%
\begin{equation}
\left\Vert f\left( t\right) \right\Vert -\func{Re}\left\langle f\left(
t\right) ,e\right\rangle \leq k\left( t\right) \ \ \ \text{for a.e. }t\in %
\left[ a,b\right] ,  \label{1.14}
\end{equation}%
then we have the inequality:%
\begin{equation}
\left( 0\leq \right) \int_{a}^{b}\left\Vert f\left( t\right) \right\Vert
dt-\left\Vert \int_{a}^{b}f\left( t\right) dt\right\Vert \leq
\int_{a}^{b}k\left( t\right) dt.  \label{1.15}
\end{equation}%
The equality holds in (\ref{1.15}) if and only if%
\begin{equation}
\int_{a}^{b}\left\Vert f\left( t\right) \right\Vert dt\geq
\int_{a}^{b}k\left( t\right) dt  \label{1.16}
\end{equation}%
and%
\begin{equation}
\int_{a}^{b}f\left( t\right) dt=\left( \int_{a}^{b}\left\Vert f\left(
t\right) \right\Vert dt-\int_{a}^{b}k\left( t\right) dt\right) e.
\label{1.17}
\end{equation}
\end{theorem}

This general result has some particular cases of interest that may be easily
applied \cite{SSD2}.

\begin{corollary}
\label{ce}If $f\in L\left( \left[ a,b\right] ;H\right) $ is such that there
exists a vector $e\in H,$ $\left\Vert e\right\Vert =1$ and $\rho \in \left(
0,1\right) $ such that 
\begin{equation}
\left\Vert f\left( t\right) -e\right\Vert \leq \rho \ \ \ \text{for a.e. }%
t\in \left[ a,b\right] ,  \label{1.18}
\end{equation}%
then 
\begin{align}
0& \leq \int_{a}^{b}\left\Vert f\left( t\right) \right\Vert dt-\left\Vert
\int_{a}^{b}f\left( t\right) dt\right\Vert   \label{1.19} \\
& \leq \frac{\rho ^{2}}{\sqrt{1-\rho ^{2}}\left( 1+\sqrt{1-\rho ^{2}}\right) 
}\func{Re}\left\langle \int_{a}^{b}f\left( t\right) dt,e\right\rangle . 
\notag
\end{align}
\end{corollary}

\begin{corollary}
\label{cf}If $f\in L\left( \left[ a,b\right] ;H\right) $ is such that there
exists a vector $e\in H,$ $\left\Vert e\right\Vert =1$ and $M\geq m>0$ such
that either%
\begin{equation}
\func{Re}\left\langle Me-f\left( t\right) ,f\left( t\right) -me\right\rangle
\geq 0\ \   \label{1.20}
\end{equation}%
or, equivalently,%
\begin{equation}
\left\Vert f\left( t\right) -\frac{M+m}{2}e\right\Vert \leq \frac{1}{2}%
\left( M-m\right) \   \label{1.21}
\end{equation}%
$\ \ $for a.e. $t\in \left[ a,b\right] ,$ then 
\begin{align}
0& \leq \int_{a}^{b}\left\Vert f\left( t\right) \right\Vert dt-\left\Vert
\int_{a}^{b}f\left( t\right) dt\right\Vert   \label{1.22} \\
& \leq \frac{\left( \sqrt{M}-\sqrt{m}\right) ^{2}}{2\sqrt{mM}}\func{Re}%
\left\langle \int_{a}^{b}f\left( t\right) dt,e\right\rangle ;  \notag
\end{align}
\end{corollary}

and finally,

\begin{corollary}
\label{cg}If $f\in L\left( \left[ a,b\right] ;H\right) $ and $r\in
L_{2}\left( \left[ a,b\right] ;H\right) ,$ $e\in H,$ $\left\Vert
e\right\Vert =1$ are such that 
\begin{equation}
\left\Vert f\left( t\right) -e\right\Vert \leq r\left( t\right) \ \ \ \text{%
for a.e. }t\in \left[ a,b\right] ,  \label{1.23}
\end{equation}%
then%
\begin{equation}
\left( 0\leq \right) \int_{a}^{b}\left\Vert f\left( t\right) \right\Vert
dt-\left\Vert \int_{a}^{b}f\left( t\right) dt\right\Vert \leq \frac{1}{2}%
\int_{a}^{b}r^{2}\left( t\right) dt.  \label{1.24}
\end{equation}
\end{corollary}

The main aim of this paper is to point out some quadratic reverses for the
continuous triangle inequality, namely, upper bounds for the nonnegative
difference%
\begin{equation*}
\left( \int_{a}^{b}\left\Vert f\left( t\right) \right\Vert dt\right)
^{2}-\left\Vert \int_{a}^{b}f\left( t\right) dt\right\Vert ^{2}
\end{equation*}%
under various assumptions on the functions $f\in L\left( \left[ a,b\right]
;H\right) .$Some related results are also pointed out. Applications for
complex-valued functions are provided as well.

\section{Quadratic Reverses of the Triangle Inequality}

The following lemma holds.

\begin{lemma}
\label{l2.1}Let $f\in L\left( \left[ a,b\right] ;H\right) $ be such that
there exists a functions $k:\Delta \subset \mathbb{R}^{2}\rightarrow \mathbb{%
R}$, $\Delta :=\left\{ \left( t,s\right) |a\leq t\leq s\leq b\right\} $ with
the property that $k\in L\left( \Delta \right) $ and%
\begin{equation}
\left( 0\leq \right) \left\Vert f\left( t\right) \right\Vert \left\Vert
f\left( s\right) \right\Vert -\func{Re}\left\langle f\left( t\right)
,f\left( s\right) \right\rangle \leq k\left( t,s\right) ,  \label{2.1}
\end{equation}%
for a.e. $\left( t,s\right) \in \Delta .$ Then we have the following
quadratic reverse of the continuous triangle inequality:%
\begin{equation}
\left( \int_{a}^{b}\left\Vert f\left( t\right) \right\Vert dt\right)
^{2}\leq \left\Vert \int_{a}^{b}f\left( t\right) dt\right\Vert
^{2}+2\iint_{\Delta }k\left( t,s\right) dtds.  \label{2.2}
\end{equation}%
The case of equality holds in (\ref{2.2}) if and only if it holds in (\ref%
{2.1}) for a.e. $\left( t,s\right) \in \Delta .$
\end{lemma}

\begin{proof}
We observe that the following identity holds%
\begin{align}
& \left( \int_{a}^{b}\left\Vert f\left( t\right) \right\Vert dt\right)
^{2}-\left\Vert \int_{a}^{b}f\left( t\right) dt\right\Vert ^{2}  \label{2.3}
\\
& =\int_{a}^{b}\int_{a}^{b}\left\Vert f\left( t\right) \right\Vert
\left\Vert f\left( s\right) \right\Vert dtds-\left\langle
\int_{a}^{b}f\left( t\right) dt,\int_{a}^{b}f\left( s\right) ds\right\rangle
\notag \\
& =\int_{a}^{b}\int_{a}^{b}\left\Vert f\left( t\right) \right\Vert
\left\Vert f\left( s\right) \right\Vert dtds-\int_{a}^{b}\int_{a}^{b}\func{Re%
}\left\langle f\left( t\right) ,f\left( s\right) \right\rangle dtds  \notag
\\
& =\int_{a}^{b}\int_{a}^{b}\left[ \left\Vert f\left( t\right) \right\Vert
\left\Vert f\left( s\right) \right\Vert -\func{Re}\left\langle f\left(
t\right) ,f\left( s\right) \right\rangle \right] dtds:=I.  \notag
\end{align}%
Now, observe that for any $\left( t,s\right) \in \left[ a,b\right] \times %
\left[ a,b\right] ,$ we have%
\begin{equation*}
\left\Vert f\left( t\right) \right\Vert \left\Vert f\left( s\right)
\right\Vert -\func{Re}\left\langle f\left( t\right) ,f\left( s\right)
\right\rangle =\left\Vert f\left( s\right) \right\Vert \left\Vert f\left(
t\right) \right\Vert -\func{Re}\left\langle f\left( s\right) ,f\left(
t\right) \right\rangle
\end{equation*}%
and thus 
\begin{equation}
I=2\iint_{\Delta }\left[ \left\Vert f\left( t\right) \right\Vert \left\Vert
f\left( s\right) \right\Vert -\func{Re}\left\langle f\left( t\right)
,f\left( s\right) \right\rangle \right] dtds.  \label{2.4}
\end{equation}%
Using the assumption (\ref{2.1}), we deduce%
\begin{equation*}
\iint_{\Delta }\left[ \left\Vert f\left( t\right) \right\Vert \left\Vert
f\left( s\right) \right\Vert -\func{Re}\left\langle f\left( t\right)
,f\left( s\right) \right\rangle \right] dtds\leq \iint_{\Delta }k\left(
t,s\right) dtds,
\end{equation*}%
and, by the identities (\ref{2.3}) and (\ref{2.4}), we deduce the desired
inequality (\ref{2.2}).

The case of equality is obvious and we omit the details.
\end{proof}

\begin{remark}
\label{r2.2}From (\ref{2.2}) one may deduce a coarser inequality that can be
useful in some applications. It is as follows:%
\begin{equation*}
\left( 0\leq \right) \int_{a}^{b}\left\Vert f\left( t\right) \right\Vert
dt-\left\Vert \int_{a}^{b}f\left( t\right) dt\right\Vert \leq \sqrt{2}\left(
\iint_{\Delta }k\left( t,s\right) dtds\right) ^{\frac{1}{2}}.
\end{equation*}
\end{remark}

\begin{remark}
\label{r2.3}If the condition (\ref{2.1}) is replaced with the following
refinement of the Schwarz inequality%
\begin{equation}
\left( 0\leq \right) k\left( t,s\right) \leq \left\Vert f\left( t\right)
\right\Vert \left\Vert f\left( s\right) \right\Vert -\func{Re}\left\langle
f\left( t\right) ,f\left( s\right) \right\rangle  \label{2.5}
\end{equation}%
for a.e. $\left( t,s\right) \in \Delta ,$ then the following refinement of
the quadratic triangle inequality is valid%
\begin{align}
\left( \int_{a}^{b}\left\Vert f\left( t\right) \right\Vert dt\right) ^{2}&
\geq \left\Vert \int_{a}^{b}f\left( t\right) dt\right\Vert
^{2}+2\iint_{\Delta }k\left( t,s\right) dtds  \label{2.6} \\
& \left( \geq \left\Vert \int_{a}^{b}f\left( t\right) dt\right\Vert
^{2}\right) .  \notag
\end{align}%
The equality holds in (\ref{2.6}) iff the case of equality holds in (\ref%
{2.5}) for a.e. $\left( t,s\right) \in \Delta .$
\end{remark}

The following result holds.

\begin{theorem}
\label{t2.4}Let $f\in L\left( \left[ a,b\right] ;H\right) $ be such that
there exists $M\geq 1\geq m\geq 0$ such that either%
\begin{equation}
\func{Re}\left\langle Mf\left( s\right) -f\left( t\right) ,f\left( t\right)
-mf\left( s\right) \right\rangle \geq 0\text{ \ for a.e. \ }\left(
t,s\right) \in \Delta ,  \label{2.7}
\end{equation}%
or, equivalently,%
\begin{equation}
\left\Vert f\left( t\right) -\frac{M+m}{2}f\left( s\right) \right\Vert \leq 
\frac{1}{2}\left( M-m\right) \left\Vert f\left( s\right) \right\Vert \text{
\ for a.e. }\left( t,s\right) \in \Delta .  \label{2.8}
\end{equation}%
Then we have the inequality:%
\begin{equation}
\left( \int_{a}^{b}\left\Vert f\left( t\right) \right\Vert dt\right)
^{2}\leq \left\Vert \int_{a}^{b}f\left( t\right) dt\right\Vert ^{2}+\frac{1}{%
2}\cdot \frac{\left( M-m\right) ^{2}}{M+m}\int_{a}^{b}\left( s-a\right)
\left\Vert f\left( s\right) \right\Vert ^{2}ds.  \label{2.9}
\end{equation}%
The case of equality holds in (\ref{2.9}) if and only if%
\begin{equation}
\left\Vert f\left( t\right) \right\Vert \left\Vert f\left( s\right)
\right\Vert -\func{Re}\left\langle f\left( t\right) ,f\left( s\right)
\right\rangle =\frac{1}{4}\cdot \frac{\left( M-m\right) ^{2}}{M+m}\left\Vert
f\left( s\right) \right\Vert ^{2}  \label{2.10}
\end{equation}%
for a.e. $\left( t,s\right) \in \Delta .$
\end{theorem}

\begin{proof}
Firstly, observe that, in an inner product space $\left( H;\left\langle
\cdot ,\cdot \right\rangle \right) $ and for $x,z,Z\in H,$ the following
statements are equivalent

\begin{enumerate}
\item[(i)] $\func{Re}\left\langle Z-x,x-z\right\rangle \geq 0$

and

\item[(ii)] $\left\Vert x-\frac{Z+z}{2}\right\Vert \leq \frac{1}{2}%
\left\Vert Z-z\right\Vert .$
\end{enumerate}

This shows that (\ref{2.7}) and (\ref{2.8}) are obviously equivalent.

Now, taking the square in (\ref{2.8}), we get%
\begin{multline*}
\left\Vert f\left( t\right) \right\Vert ^{2}+\left( \frac{M+m}{2}\right)
^{2}\left\Vert f\left( s\right) \right\Vert ^{2} \\
\leq 2\func{Re}\left\langle f\left( t\right) ,\frac{M+m}{2}f\left( s\right)
\right\rangle +\frac{1}{4}\left( M-m\right) ^{2}\left\Vert f\left( s\right)
\right\Vert ^{2},
\end{multline*}%
for a.e. $\left( t,s\right) \in \Delta ,$ and obviously, since%
\begin{equation*}
2\left( \frac{M+m}{2}\right) \left\Vert f\left( t\right) \right\Vert
\left\Vert f\left( s\right) \right\Vert \leq \left\Vert f\left( t\right)
\right\Vert ^{2}+\left( \frac{M+m}{2}\right) ^{2}\left\Vert f\left( s\right)
\right\Vert ^{2},
\end{equation*}%
we deduce that%
\begin{eqnarray*}
&&2\left( \frac{M+m}{2}\right) \left\Vert f\left( t\right) \right\Vert
\left\Vert f\left( s\right) \right\Vert  \\
&\leq &2\func{Re}\left\langle f\left( t\right) ,\frac{M+m}{2}f\left(
s\right) \right\rangle +\frac{1}{4}\left( M-m\right) ^{2}\left\Vert f\left(
s\right) \right\Vert ^{2},
\end{eqnarray*}%
giving the much simpler inequality:%
\begin{equation}
\left\Vert f\left( t\right) \right\Vert \left\Vert f\left( s\right)
\right\Vert -\func{Re}\left\langle f\left( t\right) ,f\left( s\right)
\right\rangle \leq \frac{1}{4}\cdot \frac{\left( M-m\right) ^{2}}{M+m}%
\left\Vert f\left( s\right) \right\Vert ^{2}  \label{2.11}
\end{equation}%
for a.e. $\left( t,s\right) \in \Delta .$

Applying Lemma \ref{l2.1} for $k\left( t,s\right) :=\frac{1}{4}\cdot \frac{%
\left( M-m\right) ^{2}}{M+m}\left\Vert f\left( s\right) \right\Vert ^{2},$
we deduce%
\begin{equation}
\left( \int_{a}^{b}\left\Vert f\left( t\right) \right\Vert dt\right)
^{2}\leq \left\Vert \int_{a}^{b}f\left( t\right) dt\right\Vert ^{2}+\frac{1}{%
2}\cdot \frac{\left( M-m\right) ^{2}}{M+m}\iint_{\Delta }\left\Vert f\left(
s\right) \right\Vert ^{2}ds  \label{2.12}
\end{equation}%
with equality if and only if (\ref{2.11}) holds for a.e. $\left( t,s\right)
\in \Delta .$

Since%
\begin{equation*}
\iint_{\Delta }\left\Vert f\left( s\right) \right\Vert
^{2}ds=\int_{a}^{b}\left( \int_{a}^{s}\left\Vert f\left( s\right)
\right\Vert ^{2}dt\right) ds=\int_{a}^{b}\left( s-a\right) \left\Vert
f\left( s\right) \right\Vert ^{2}ds,
\end{equation*}%
then by (\ref{2.12}) we deduce the desired result (\ref{2.9}).
\end{proof}

Another result which is similar to the one above is incorporated in the
following theorem.

\begin{theorem}
\label{t2.5}With the assumptions of Theorem \ref{t2.4}, we have%
\begin{equation}
\left( \int_{a}^{b}\left\Vert f\left( t\right) \right\Vert dt\right)
^{2}-\left\Vert \int_{a}^{b}f\left( t\right) dt\right\Vert ^{2}\leq \frac{%
\left( \sqrt{M}-\sqrt{m}\right) ^{2}}{2\sqrt{Mm}}\left\Vert
\int_{a}^{b}f\left( t\right) dt\right\Vert ^{2}  \label{2.13}
\end{equation}%
or, equivalently,%
\begin{equation}
\int_{a}^{b}\left\Vert f\left( t\right) \right\Vert dt\leq \left( \frac{M+m}{%
2\sqrt{Mm}}\right) ^{\frac{1}{2}}\left\Vert \int_{a}^{b}f\left( t\right)
dt\right\Vert .  \label{2.14}
\end{equation}%
The case of equality holds in (\ref{2.13}) or (\ref{2.14}) if and only if 
\begin{equation}
\left\Vert f\left( t\right) \right\Vert \left\Vert f\left( s\right)
\right\Vert =\frac{M+m}{2\sqrt{Mm}}\func{Re}\left\langle f\left( t\right)
,f\left( s\right) \right\rangle ,  \label{2.14'}
\end{equation}%
for a.e. $\left( t,s\right) \in \Delta .$
\end{theorem}

\begin{proof}
From (\ref{2.7}), we deduce%
\begin{equation*}
\left\Vert f\left( t\right) \right\Vert ^{2}+Mm\left\Vert f\left( s\right)
\right\Vert ^{2}\leq \left( M+m\right) \func{Re}\left\langle f\left(
t\right) ,f\left( s\right) \right\rangle 
\end{equation*}%
for a.e. $\left( t,s\right) \in \Delta .$ Dividing by $\sqrt{Mm}>0,$ we
deduce%
\begin{equation*}
\frac{\left\Vert f\left( t\right) \right\Vert ^{2}}{\sqrt{Mm}}+\sqrt{Mm}%
\left\Vert f\left( s\right) \right\Vert ^{2}\leq \frac{M+m}{\sqrt{Mm}}\func{%
Re}\left\langle f\left( t\right) ,f\left( s\right) \right\rangle 
\end{equation*}%
and, obviously, since%
\begin{equation*}
2\left\Vert f\left( t\right) \right\Vert \left\Vert f\left( s\right)
\right\Vert \leq \frac{\left\Vert f\left( t\right) \right\Vert ^{2}}{\sqrt{Mm%
}}+\sqrt{Mm}\left\Vert f\left( s\right) \right\Vert ^{2},
\end{equation*}%
hence%
\begin{equation*}
\left\Vert f\left( t\right) \right\Vert \left\Vert f\left( s\right)
\right\Vert \leq \frac{M+m}{\sqrt{Mm}}\func{Re}\left\langle f\left( t\right)
,f\left( s\right) \right\rangle 
\end{equation*}%
for a.e. $\left( t,s\right) \in \Delta ,$ giving%
\begin{equation*}
\left\Vert f\left( t\right) \right\Vert \left\Vert f\left( s\right)
\right\Vert -\func{Re}\left\langle f\left( t\right) ,f\left( s\right)
\right\rangle \leq \frac{\left( \sqrt{M}-\sqrt{m}\right) ^{2}}{2\sqrt{Mm}}%
\func{Re}\left\langle f\left( t\right) ,f\left( s\right) \right\rangle .
\end{equation*}%
Applying Lemma \ref{l2.1} for $k\left( t,s\right) :=\frac{\left( \sqrt{M}-%
\sqrt{m}\right) ^{2}}{\sqrt{Mm}}\func{Re}\left\langle f\left( t\right)
,f\left( s\right) \right\rangle ,$ we deduce%
\begin{equation}
\left( \int_{a}^{b}\left\Vert f\left( t\right) \right\Vert dt\right)
^{2}\leq \left\Vert \int_{a}^{b}f\left( t\right) dt\right\Vert ^{2}+\frac{%
\left( \sqrt{M}-\sqrt{m}\right) ^{2}}{2\sqrt{Mm}}\func{Re}\left\langle
f\left( t\right) ,f\left( s\right) \right\rangle .  \label{2.15}
\end{equation}%
On the other hand, since%
\begin{equation*}
\func{Re}\left\langle f\left( t\right) ,f\left( s\right) \right\rangle =%
\func{Re}\left\langle f\left( s\right) ,f\left( t\right) \right\rangle \text{
\ for any \ }\left( t,s\right) \in \left[ a,b\right] ^{2},
\end{equation*}%
hence%
\begin{align*}
\iint_{\Delta }\func{Re}\left\langle f\left( t\right) ,f\left( s\right)
\right\rangle dtds& =\frac{1}{2}\int_{a}^{b}\int_{a}^{b}\func{Re}%
\left\langle f\left( t\right) ,f\left( s\right) \right\rangle dtds \\
& =\frac{1}{2}\func{Re}\left\langle \int_{a}^{b}f\left( t\right)
dt,\int_{a}^{b}f\left( s\right) ds\right\rangle  \\
& =\frac{1}{2}\left\Vert \int_{a}^{b}f\left( t\right) dt\right\Vert ^{2}
\end{align*}%
and thus, from (\ref{2.15}), we get (\ref{2.13}).

The equivalence between (\ref{2.13}) and (\ref{2.14}) is obvious and we omit
the details.
\end{proof}

\section{Related Results}

The following result also holds.

\begin{theorem}
\label{t3.1}Let $f\in L\left( \left[ a,b\right] ;H\right) $ and $\gamma
,\Gamma \in \mathbb{R}$ be such that either%
\begin{equation}
\func{Re}\left\langle \Gamma f\left( s\right) -f\left( t\right) ,f\left(
t\right) -\gamma f\left( s\right) \right\rangle \geq 0\text{ \ for a.e. \ }%
\left( t,s\right) \in \Delta ,  \label{3.1}
\end{equation}%
or, equivalently,%
\begin{equation}
\left\Vert f\left( t\right) -\frac{\Gamma +\gamma }{2}f\left( s\right)
\right\Vert \leq \frac{1}{2}\left\vert \Gamma -\gamma \right\vert \left\Vert
f\left( s\right) \right\Vert \text{ \ for a.e. \ }\left( t,s\right) \in
\Delta .  \label{3.2}
\end{equation}%
Then we have the inequality:%
\begin{equation}
\int_{a}^{b}\left[ \left( b-s\right) +\gamma \Gamma \left( s-a\right) \right]
\left\Vert f\left( s\right) \right\Vert ^{2}ds\leq \frac{\Gamma +\gamma }{2}%
\left\Vert \int_{a}^{b}f\left( s\right) ds\right\Vert ^{2}.  \label{3.3}
\end{equation}%
The case of equality holds in (\ref{3.3}) if and only if the case of
equality holds in either (\ref{3.1}) or (\ref{3.2}) for a.e. $\left(
t,s\right) \in \Delta $.
\end{theorem}

\begin{proof}
The inequality (\ref{3.1}) is obviously equivalent to%
\begin{equation}
\left\Vert f\left( t\right) \right\Vert ^{2}+\gamma \Gamma \left\Vert
f\left( s\right) \right\Vert ^{2}\leq \left( \Gamma +\gamma \right) \func{Re}%
\left\langle f\left( t\right) ,f\left( s\right) \right\rangle  \label{3.4}
\end{equation}%
for a.e. $\left( t,s\right) \in \Delta .$

Integrating (\ref{3.4}) on $\Delta ,$ we deduce%
\begin{multline}
\int_{a}^{b}\left( \int_{a}^{s}\left\Vert f\left( t\right) \right\Vert
^{2}dt\right) ds+\gamma \Gamma \int_{a}^{b}\left( \left\Vert f\left(
s\right) \right\Vert ^{2}\int_{a}^{s}dt\right) ds  \label{3.5} \\
=\left( \Gamma +\gamma \right) \int_{a}^{b}\left( \int_{a}^{s}\func{Re}%
\left\langle f\left( t\right) ,f\left( s\right) \right\rangle dt\right) ds.
\end{multline}%
It is easy to see, on integrating by parts, that%
\begin{align*}
\int_{a}^{b}\left( \int_{a}^{s}\left\Vert f\left( t\right) \right\Vert
^{2}dt\right) ds& =s\left. \int_{a}^{s}\left\Vert f\left( t\right)
\right\Vert ^{2}dt\right\vert _{a}^{b}-\int_{a}^{b}s\left\Vert f\left(
s\right) \right\Vert ^{2}ds \\
& =b\int_{a}^{s}\left\Vert f\left( s\right) \right\Vert
^{2}ds-\int_{a}^{b}s\left\Vert f\left( s\right) \right\Vert ^{2}ds \\
& =\int_{a}^{b}\left( b-s\right) \left\Vert f\left( s\right) \right\Vert
^{2}ds
\end{align*}%
and 
\begin{equation*}
\int_{a}^{b}\left( \left\Vert f\left( s\right) \right\Vert
^{2}\int_{a}^{s}dt\right) ds=\int_{a}^{b}\left( s-a\right) \left\Vert
f\left( s\right) \right\Vert ^{2}ds.
\end{equation*}%
Since%
\begin{align*}
\frac{d}{ds}\left( \left\Vert \int_{a}^{b}f\left( t\right) dt\right\Vert
^{2}\right) & =\frac{d}{ds}\left\langle \int_{a}^{s}f\left( t\right)
dt,\int_{a}^{s}f\left( t\right) dt\right\rangle  \\
& =\left\langle f\left( s\right) ,\int_{a}^{s}f\left( t\right)
dt\right\rangle +\left\langle \int_{a}^{s}f\left( t\right) dt,f\left(
s\right) \right\rangle  \\
& =2\func{Re}\left\langle \int_{a}^{s}f\left( t\right) dt,f\left( s\right)
\right\rangle ,
\end{align*}%
hence%
\begin{align*}
\int_{a}^{b}\left( \int_{a}^{s}\func{Re}\left\langle f\left( t\right)
,f\left( s\right) \right\rangle dt\right) ds& =\int_{a}^{b}\func{Re}%
\left\langle \int_{a}^{s}f\left( t\right) dt,f\left( s\right) \right\rangle
ds \\
& =\frac{1}{2}\int_{a}^{b}\frac{d}{ds}\left( \left\Vert \int_{a}^{s}f\left(
t\right) dt\right\Vert ^{2}\right) ds \\
& =\frac{1}{2}\left. \left\Vert \int_{a}^{s}f\left( t\right) dt\right\Vert
^{2}\right\vert _{a}^{b} \\
& =\frac{1}{2}\left\Vert \int_{a}^{b}f\left( t\right) dt\right\Vert ^{2}.
\end{align*}%
Utilising (\ref{3.5}), we deduce the desired inequality (\ref{3.3}).

The case of equality is obvious and we omit the details.
\end{proof}

\begin{remark}
Consider the function $\varphi \left( s\right) :=\left( b-s\right) +\gamma
\Gamma \left( s-a\right) ,$ $s\in \left[ a,b\right] .$ Obviously,%
\begin{equation*}
\varphi \left( s\right) =\left( \Gamma \gamma -1\right) s+b-\gamma \Gamma a.
\end{equation*}%
Observe that, if $\Gamma \gamma \geq 1,$ then%
\begin{equation*}
b-a=\varphi \left( a\right) \leq \varphi \left( s\right) \leq \varphi \left(
b\right) =\gamma \Gamma \left( b-a\right) ,\ \ \ \ s\in \left[ a,b\right]
\end{equation*}%
and, if $\Gamma \gamma <1,$ then%
\begin{equation*}
\gamma \Gamma \left( b-a\right) \leq \varphi \left( s\right) \leq b-a,\ \ \
\ s\in \left[ a,b\right] .
\end{equation*}
\end{remark}

Taking into account the above remark, we may state the following corollary.

\begin{corollary}
\label{c3.2}Assume that $f,\gamma ,\Gamma $ are as in Theorem \ref{t3.1}.

\begin{enumerate}
\item[a)] If $\Gamma \gamma \geq 1,$ then we have the inequality%
\begin{equation*}
\left( b-a\right) \int_{a}^{b}\left\Vert f\left( s\right) \right\Vert
^{2}ds\leq \frac{\Gamma +\gamma }{2}\left\Vert \int_{a}^{b}f\left( s\right)
ds\right\Vert ^{2}.
\end{equation*}

\item[b)] If $0<\Gamma \gamma <1,$ then we have the inequality%
\begin{equation*}
\gamma \Gamma \left( b-a\right) \int_{a}^{b}\left\Vert f\left( s\right)
\right\Vert ^{2}ds\leq \frac{\Gamma +\gamma }{2}\left\Vert
\int_{a}^{b}f\left( s\right) ds\right\Vert ^{2}.
\end{equation*}
\end{enumerate}
\end{corollary}

\section{Applications for Complex-Valued Functions}

Let $f:\left[ a,b\right] \rightarrow \mathbb{C}$ be a Lebesgue integrable
function and $M\geq 1\geq m\geq 0.$ The condition (\ref{2.7}) from Theorem %
\ref{t2.4}, which plays a fundamental role in the results obtained above,
can be translated in this case as%
\begin{equation}
\func{Re}\left[ \left( Mf\left( s\right) -f\left( t\right) \right) \left( 
\overline{f\left( t\right) }-m\overline{f\left( s\right) }\right) \right]
\geq 0  \label{e.4.1}
\end{equation}%
for a.e. $a\leq t\leq s\leq b.$

Since, obviously%
\begin{eqnarray*}
&&\func{Re}\left[ \left( Mf\left( s\right) -f\left( t\right) \right) \left( 
\overline{f\left( t\right) }-m\overline{f\left( s\right) }\right) \right]  \\
&=&\left[ \left( M\func{Re}f\left( s\right) -\func{Re}f\left( t\right)
\right) \left( \func{Re}f\left( t\right) -m\func{Re}f\left( s\right) \right) %
\right]  \\
&&+\left[ \left( M\func{Im}f\left( s\right) -\func{Im}f\left( t\right)
\right) \left( \func{Im}f\left( t\right) -m\func{Im}f\left( s\right) \right) %
\right] 
\end{eqnarray*}%
hence a sufficient condition for the inequality in (\ref{e.4.1}) to hold is%
\begin{equation}
m\func{Re}f\left( s\right) \leq \func{Re}f\left( t\right) \leq M\func{Re}%
f\left( s\right) \text{ and }m\func{Im}f\left( s\right) \leq \func{Im}%
f\left( t\right) \leq M\func{Im}f\left( s\right)   \label{e.4.2}
\end{equation}%
for a.e. $a\leq t\leq s\leq b.$

Utilising Theorems \ref{t2.4},\ref{t2.5} and \ref{t3.1} we may state the
following results incorporating quadratic reverses of the continuous
triangle inequality:

\begin{proposition}
\label{p.4.1} With the above assumptions for $f,M$ and $m,$ and if (\ref%
{e.4.1}) holds true, then we have the inequalities%
\begin{equation*}
\left( \int_{a}^{b}\left\vert f\left( t\right) \right\vert dt\right)
^{2}\leq \left\vert \int_{a}^{b}f\left( t\right) dt\right\vert ^{2}+\frac{1}{%
2}\cdot \frac{\left( M-m\right) ^{2}}{M+m}\int_{a}^{b}\left( s-a\right)
\left\vert f\left( s\right) \right\vert ^{2}ds,
\end{equation*}%
\begin{equation*}
\int_{a}^{b}\left\vert f\left( t\right) \right\vert dt\leq \left( \frac{M+m}{%
2\sqrt{Mm}}\right) ^{\frac{1}{2}}\left\vert \int_{a}^{b}f\left( t\right)
dt\right\vert ,
\end{equation*}%
and%
\begin{equation*}
\int_{a}^{b}\left[ \left( b-s\right) +\gamma \Gamma \left( s-a\right) \right]
\left\vert f\left( s\right) \right\vert ^{2}ds\leq \frac{\Gamma +\gamma }{2}%
\left\vert \int_{a}^{b}f\left( s\right) ds\right\vert ^{2}.
\end{equation*}
\end{proposition}

\begin{remark}
\label{r.4.1} One may wonder if there are functions satisfying the condition
(\ref{e.4.2}) above. It suffices to find examples of real  functions $%
\varphi :\left[ a,b\right] \rightarrow \mathbb{R}$ verifying the following
double inequality%
\begin{equation}
\gamma \varphi \left( s\right) \leq \varphi \left( t\right) \leq \Gamma
\varphi \left( s\right)   \label{e.4.3}
\end{equation}%
for some given $\gamma ,\Gamma $ with $0\leq \gamma \leq 1\leq \Gamma
<\infty $ for a.e. $a\leq t\leq s\leq b.$

For this purpose, consider $\psi :\left[ a,b\right] \rightarrow \mathbb{R}$
a differentiable function on $\left( a,b\right) $, continuous on $\left[ a,b%
\right] $ and with the property that there exists $\Theta \geq 0\geq \theta $
such that%
\begin{equation}
\theta \leq \psi ^{\prime }\left( u\right) \leq \Theta \text{ for any }u\in
\left( a,b\right) .  \label{e.4.4}
\end{equation}%
By Lagrange's mean value theorem, we have, for any $a\leq t\leq s\leq b$%
\begin{equation*}
\psi \left( s\right) -\psi \left( t\right) =\psi ^{\prime }\left( \xi
\right) \left( s-t\right) 
\end{equation*}%
with $t\leq \xi \leq s.$ Therefore, for $a\leq t\leq s\leq b,$ by (\ref%
{e.4.4}), we have the inequality%
\begin{equation*}
\theta \left( b-a\right) \leq \theta \left( s-t\right) \leq \psi \left(
s\right) -\psi \left( t\right) \leq \Theta \left( s-t\right) \leq \Theta
\left( b-a\right) .
\end{equation*}%
If we choose the function $\varphi :\left[ a,b\right] \rightarrow \mathbb{R}$
given by%
\begin{equation*}
\varphi \left( t\right) :=\exp \left[ -\psi \left( t\right) \right] ,\text{ }%
t\in \left[ a,b\right] ,
\end{equation*}%
and $\gamma :=\exp \left[ \theta \left( b-a\right) \right] \leq 1,$ $\Gamma
:=\exp \left[ \Theta \left( b-a\right) \right] ,$ then (\ref{e.4.3}) holds
true for any $a\leq t\leq s\leq b.$
\end{remark}

\end{document}